\documentclass[11pt]{amsart}

\usepackage{amsmath}
\usepackage{amssymb}
\usepackage{mathrsfs}
\usepackage[all]{xy}
\usepackage{comment}
\usepackage{color}
\usepackage[colorlinks,citecolor=blue,urlcolor=black,linkcolor=black]{hyperref}

\newtheorem{theorem}{Theorem}[section]

\newtheorem{corollary}[theorem]{Corollary}

\theoremstyle{definition}

\theoremstyle{remark}

\newtheorem{conjecture}[theorem]{Conjecture}

\newcount\skewfactor
\def\mathunderaccent#1#2 {\let\theaccent#1\skewfactor#2
\mathpalette\putaccentunder}
\def\putaccentunder#1#2{\oalign{$#1#2$\crcr\hidewidth
\vbox to.2ex{\hbox{$#1\skew\skewfactor\theaccent{}$}\vss}\hidewidth}}

\def\smallbox#1{\leavevmode\thinspace\hbox{\vrule\vtop{\vbox
   {\hrule\kern1pt\hbox{\vphantom{\tt/}\thinspace{\tt#1}\thinspace}}
   \kern1pt\hrule}\vrule}\thinspace}

\newcommand{\cf}{{\rm cf}}

\newcommand{\stick}{\ensuremath\mspace{2mu}\mid \mspace{-12mu}{\raise 0.4em \hbox{$\bullet$}}}

\def\qedref#1{$\qed_{\reforiginal{#1}}$}

\title{Weak diamond and pcf theory}
\author{Shimon Garti}
\address{Einstein Institute of Mathematics,
 The Hebrew University of Jerusalem,
 Jerusalem 91904, Israel}
\email{shimon.garty@mail.huji.ac.il}
\thanks{}


\begin{document}
\let\labeloriginal\label
\let\reforiginal\ref

\begin{abstract}
We obtain bounds on the cardinality of $pcf(\mathfrak{a})$ from instances of weak diamond.
Consequently, under mild assumptions there are many singular cardinals of the form $\aleph_\delta$ for which $2^{\aleph_\delta}<\aleph_{|\delta|^{+3}}$.
For example, if every limit cardinal is a strong limit cardinal then this bound holds at a class of singular cardinals.
\end{abstract}

\maketitle

\newpage

\section{Introduction}

The singular cardinals problem is the question of possible values of $2^\lambda$, where $\lambda$ is a strong limit singular cardinal.
The problem focuses on singular cardinals since the behaviour of $2^\kappa$ when $\kappa$ is regular is very well understood.
In fact, the global demeanor of the power set operation at every regular cardinal is fully described by Easton's theorem from \cite{MR0269497}.

The assumption of strong limitude comes from the fact that one can increase $2^\lambda$ (when $\lambda$ is singular) simply by choosing a regular cardinal $\theta\in\lambda$ and increasing $2^\theta$.
In order to avoid this artificial influence on $2^\lambda$ one has to assume that no cardinal below $\lambda$ forces a value for $2^\lambda$, that is, $2^\theta\in\lambda$ whenever $\theta\in\lambda$.
In other words, $\lambda$ is a strong limit cardinal.

One of the highlights of set theory in the twentieth century is the following theorem of Shelah from \cite{MR1318912}.
If $\aleph_\delta$ is a strong limit cardinal and $\delta<\aleph_\delta$ then $2^{\aleph_\delta}<\aleph_{|\delta|^{+4}}$.
In order to prove this theorem, Shelah developed pcf theory and obtained many results in cardinal arithmetic, including this bound on the power set of strong limit singular cardinals.

To place these results in context we should mention here the other side of the coin.
Suppose, for simplicity, that $\lambda>\cf(\lambda)=\kappa$ and $2^\beta=\beta^+$ for every $\beta\in\lambda$.
If $\kappa>\aleph_0$ then $2^\lambda=\lambda^+$, as proved by Silver in \cite{MR0429564}.\footnote{A combinatorial proof of Silver's theorem appeared in \cite{MR0387061}. The analysis of $2^\lambda$ in that paper resembles basic features of pcf theory.}
But Silver's argument is not applicable if $\kappa=\aleph_0$.
In fact, not only the argument but also the statement.
Magidor proved in \cite{MR491183}, \cite{MR491184} that consistently $2^{\aleph_0}=\aleph_{n+1}$ for every $n\in\omega$ but $2^{\aleph_\omega}>\aleph_{\omega+1}$.
It should be emphasized that the value of $2^{\aleph_\omega}$ above $\aleph_{\omega+1}$ obtained by Magidor was relatively small.
In fact, $2^{\aleph_\omega}<\aleph_{\omega_1}$ in all known models of the failure of \textsf{SCH} at $\aleph_\omega$,\footnote{\textsf{SCH} is an acronym for the singular cardinals hypothesis. A common formulation is $2^\lambda=\lambda^+$ whenever $\lambda$ is a strong limit singular cardinal.} and the same is true at every strong limit singular cardinal that is not a fixed point of the aleph function.
Shelah's theorem in the specific case of $\aleph_\omega$ amounts to $2^{\aleph_\omega}<\aleph_{\omega_4}$, so we have an interesting gap between the upper bound and the ability to increase $2^{\aleph_\omega}$ under the above assumption.

Despite a lot of effort in the last few decades, no significant progress was made in either direction.
Our goal in this paper is to consider an improvement of Shelah's bound.
That is, we try to reduce $|\delta|^{+4}$ to $|\delta|^{+3}$.
The idea is that if the class of all strong limit singular cardinals is deemed then many of them satisfy the inequality $2^{\aleph_\delta}<\aleph_{|\delta|^{+3}}$, where many of them means class-many.

Our results hinge upon some cardinal arithmetic assumptions, but these assumptions are relatively weak, and we actually believe that the main result is provable in \textsf{ZFC}.
Here is a typical statement which we can prove.
Suppose that every limit cardinal is a strong limit cardinal.
Then there is a class of singular strong limit cardinals of the form $\aleph_\delta$, such that $2^{\aleph_\delta}<\aleph_{|\delta|^{+3}}$.
We emphasize that this conclusion follows from much a weaker assumption, as will be explicated later.

The paper contains three additional sections.
In the first one we give some background, focusing on weak diamond and on pcf theory.
In the second we describe Galvin's property and then we state and prove a club guessing theorem which follows from Galvin's property.
In the last one we derive the main results in cardinal arithmetic.

Our notation is (hopefully) standard.
If $\kappa=\cf(\kappa)<\lambda$ then $S^\lambda_\kappa=\{\delta\in\lambda\mid\cf(\delta)=\kappa\}$.
This set is a stationary subset of $\lambda$, provided that $\cf(\lambda)>\omega$.
If $E$ is a subset of $\kappa$ then $acc(E)=\{\delta\in\kappa\mid \bigcup(E\cap\delta)=\delta\}$.
This set is called \emph{the set of accumulation points} of $E$, and in the main application $E$ is a club of $\kappa$, in which case $acc(E)$ is also a club of $\kappa$.
Other notation will be introduced as the need arises.

\newpage

\section{Preliminaries}

In this section we touch upon two topics.
The first is a prediction principle belonging to the diamond family, and the second is pcf theory.
Let us commence with the Devlin-Shelah weak diamond from \cite{MR0469756}.

Recall that a diamond sequence (at $\aleph_1$) is a sequence of sets $(A_\alpha\mid\alpha\in\omega_1)$ such that for every $A\subseteq\omega_1$ the set $S_A=\{\alpha\in\omega_1\mid A\cap\alpha=A_\alpha\}$ is a stationary subset of $\aleph_1$.
The diamond principle $\Diamond_{\aleph_1}$, discovered by Jensen in \cite{MR0309729}, is the statement that there exists a diamond sequence.
It is easy to see that $\Diamond_{\aleph_1}$ implies $2^{\aleph_0}=\aleph_1$.
A deep ancient result of Jensen shows that $\Diamond_{\aleph_1}$ is strictly stronger than $2^{\aleph_0}=\aleph_1$, see \cite{MR0384542} for a detailed account of a forcing construction of this result, namely the failure of the principle $\Diamond_{\aleph_1}$ in a model of \textsf{CH}.

Motivated by algebraic problems, Devlin and Shelah phrased a prediction principle that is sufficiently strong to imply some of the consequences of diamond but sufficiently weak to follow from \textsf{CH}.
The weak diamond $\Phi_{\aleph_1}$ is the statement that for every $c:{}^{<\omega_1}2\rightarrow{2}$ one can find $g\in{}^{\omega_1}2$ so that for every $f\in{}^{\omega_1}2$ the set $\{\alpha\in\omega_1\mid c(f\upharpoonright\alpha)=g(\alpha)\}$ is stationary in $\omega_1$.
It follows that $\Phi_{\aleph_1}$ implies $2^{\aleph_0}<2^{\aleph_1}$ and vice versa.
So, unlike the diamond, here the prediction principle at $\aleph_1$ is equivalent to the cardinal arithmetic statement.\footnote{The statement $2^{\aleph_0}<2^{\aleph_1}$ is called sometimes \emph{the weak continuum hypothesis}.}

This equivalence holds true in general.
One can define $\Phi_{\kappa}$ upon replacing $\aleph_1$ by $\kappa$ in the definition of the weak diamond, and then $2^\kappa<2^{\kappa^+}$ is equivalent to $\Phi_{\kappa^+}$ at every infinite cardinal $\kappa$.\footnote{The substantial direction can be deduced from \cite{MR0469756} by replacing any occurrence of $\aleph_0$ with $\kappa$ and any occurrence of $\aleph_1$ with $\kappa^+$. The proof of the easy direction is spelled-out in \cite{MR3604115}.}
Modern research shows that for $\kappa>\aleph_0$, the local instance of \textsf{GCH} expressed by $2^\kappa=\kappa^+$ is equivalent to $\Diamond_{\kappa^+}$, so diamond and weak diamond become similar from this point of view.

There is, however, a crucial difference between these principles if one adopts a global point of view.
It is consistent that \textsf{GCH} fails everywhere, and then of course $\Diamond_{\kappa^+}$ fails for every infinite cardinal $\kappa$.\footnote{Let us indicate that $\Diamond_\kappa$ holds if $\kappa$ is a sufficiently large cardinal, e.g. if $\kappa$ is measurable. But if one concentrates on successor cardinals then $\Diamond_{\kappa^+}$ fails everywhere in models of global failure of \textsf{GCH}, e.g. \cite{MR1087344}.}
But weak diamond most hold, in \textsf{ZFC}, at a class of infinite cardinals.
The following theorem appears as \cite[Proposition 2.14]{MR4611828} with a short sketch of the proof.
Let us state the theorem and give full details.

\begin{theorem}
  \label{thmclassofwd} Let $\kappa_0$ be an infinite cardinal.
  Then the weak diamond $\Phi_\kappa$ holds at some $\kappa>\kappa_0$.
\end{theorem}

\par\noindent\emph{Proof}. \newline
Let $\mu=2^{\kappa_0}$.
We consider three possible cases.
In the first case, $\mu$ is a successor cardinal, say $\mu=\theta^+$.
Now if $2^\theta=\theta^+=\mu$ then $2^\theta<2^{\theta^+}$ since $2^{\theta^+}=2^\mu>\mu$, and then $\Phi_{\theta^+}$ holds by \cite{MR0469756}.

If $2^\theta>\mu$ then let $\kappa$ be the first cardinal for which $2^\kappa>\mu$, and notice that $\kappa>\kappa_0$.
Now if $\kappa$ is a successor cardinal then (letting $\kappa=\chi^+$) one has $2^\chi<2^\kappa=2^{\chi^+}$ and $\Phi_{\chi^+}$ holds.
If $\kappa$ is a limit cardinal then, since $\kappa$ is the first cardinal which satisfies $2^\kappa>\mu$ and since $2^{\kappa_0}=\mu$, we see that $2^\tau=\mu$ for every $\tau\in[\kappa_0,\kappa)$.
It follows that $\kappa$ cannot be singular,\footnote{This is a consequence of the Bukovsk\'y-Hechler theorem, see \cite{MR0183649}.} thus $\kappa$ is a regular limit cardinal.
However, $\kappa_0<\kappa$ and $2^{\kappa_0}=\mu>\kappa$, hence $\kappa$ is weakly and not strongly inaccessible cardinal.
In this case, $\Phi_\kappa$ holds\footnote{The statement appears without proof in \cite{MR0469756}. The detailed argument can be found in \cite[Theorem 1.3]{MR3914938}.}.
We conclude, therefore, that in all subcases of the case in which $\mu=2^{\kappa_0}$ is a successor cardinal one has weak diamond at some $\kappa>\kappa_0$, so the first case is covered.

In the second case, $\mu=2^{\kappa_0}$ is a singular cardinal.
Let $\theta=\cf(\mu)$.
Necessarily, $\theta>\kappa_0$ and hence $2^\theta\geq\mu$.
But $\cf(2^\theta)>\theta$, so $2^\theta>\mu$.
Let $\theta_0$ be the first cardinal so that $2^{\theta_0}>\mu$.
Either $\theta_0$ is a successor cardinal or weakly inaccessible.
In both cases $\Phi_{\theta_0}$ holds, as explained in the previous case.
Since $\theta_0>\kappa_0$ (recall that $2^{\kappa_0}=\mu<2^{\theta_0}$), the second case is covered as well.

The last possible case is when $\mu=2^{\kappa_0}$ is weakly (but not strongly) inaccessible.
If $2^\theta=\mu$ for every $\theta\in[\kappa_0,\mu)$ then $\Phi_\mu$ holds by \cite{MR0469756}.
If not, let $\kappa$ be the first cardinal for which $2^\kappa>\mu$.
Notice that $\kappa>\kappa_0$.
As in the previous cases, either $\kappa$ is a successor cardinal or weakly inaccessible.
In both alternatives, $\Phi_\kappa$ holds as mentioned before, so we are done.

\hfill \qedref{thmclassofwd}

In the rest of this section we survey some basic facts and definitions from pcf theory.
Motivated by the singular cardinals problem, Shelah developed pcf theory in order to determine the possible values of $2^\lambda$ when $\lambda$ is a strong limit singular cardinal.
Very quickly, Shelah realized that one should understand the possible cofinalities of products of regular cardinals below $\lambda$, in many cases end-segments of ${\rm Reg}\cap\lambda$ (where ${\rm Reg}$ stands for the class of regular cardinals).

Let $\mathfrak{a}$ be a set of regular cardinals.
In most theorems one has to assume that $|\mathfrak{a}|<\min(\mathfrak{a})$.
A set of regular cardinals which satisfies this proviso will be called \emph{progressive}.
If $\lambda$ is a fixed point of the aleph-function, that is $\lambda=\aleph_\delta=\delta$, then every end-segment $\mathfrak{a}$ of ${\rm Reg}\cap\lambda$ is of size $\lambda$.
Thus, we usually assume that $\lambda$ is not a fixed point of the aleph-function.

Let $\mathfrak{a}=\{\lambda_i\mid i\in\kappa\}$ be a progressive set.
Let $\mathcal{J}$ be an ideal over $\kappa$.\footnote{We always make the assumption that $\mathcal{J}\supseteq\mathcal{J}^{\rm bd}_\kappa$, where $\mathcal{J}^{\rm bd}_\kappa$ denotes the ideal of bounded subsets of $\kappa$.}
The ideal $\mathcal{J}$ gives rise to a partial ordering defined on $\prod\mathfrak{a}$ as follows.
For $f,g\in\prod\mathfrak{a}$ one says that $f<_{\mathcal{J}}g$ iff $\{i\in\kappa\mid f(i)\geq g(i)\}\in\mathcal{J}$.
A sequence $\bar{f}=(f_\alpha\mid\alpha\in\lambda)$ of elements of $\prod\mathfrak{a}$ is \emph{a scale} in $(\prod\mathfrak{a},\mathcal{J})$ iff $\bar{f}$ is both \emph{increasing} (to wit, $\alpha<\beta<\lambda\Rightarrow f_\alpha<_\mathcal{J}f_\beta$) and \emph{cofinal} (that is, for every $h\in\prod\mathfrak{a}$ there is $\alpha\in\lambda$ such that $h<_{\mathcal{J}}f_\alpha$).

The product $(\prod\mathfrak{a},\mathcal{J})$ has \emph{true cofinality} iff there is a scale in $(\prod\mathfrak{a},\mathcal{J})$, in which case the true cofinality of $(\prod\mathfrak{a},\mathcal{J})$ is the minimal length of such a scale.
We shall write $tcf(\prod\mathfrak{a},\mathcal{J})=\lambda$.
Given a progressive set $\mathfrak{a}$, $pcf(\mathfrak{a})$ is the set of all $tcf(\prod\mathfrak{a},\mathcal{J})$ where $\mathcal{J}$ is an ideal over $\mathfrak{a}$.
It follows immediately that $\mathfrak{a}\subseteq pcf(\mathfrak{a})$ and hence $|pcf(\mathfrak{a})|\geq|\mathfrak{a}|$.
A central challenge in pcf theory is to find an upper bound on the cardinality of $pcf(\mathfrak{a})$, especially when $\mathfrak{a}$ is an interval of regular cardinals.

By finding such an upper bound, one can derive non-trivial conclusions with regard to the singular cardinals problem.
The reason is that $pcf(\mathfrak{a})$ has always a last element dubbed as $\max pcf(\mathfrak{a})$.
A fundamental theorem of Shelah says that if $\lambda$ is a strong limit singular cardinal that is not a fixed point of the $\aleph$-function, and if $\mathfrak{a}$ is a progressive end-segment of ${\rm Reg}\cap\lambda$, then $2^\lambda=\max pcf(\mathfrak{a})$.
Therefore, an upper bound on $pcf(\mathfrak{a})$ would give an upper bound on the value of $2^\lambda$.
It should be noted that an upper bound in terms of the $\beth$-function also exists,\footnote{Such bounds on $2^\lambda$ when $\lambda>\cf(\lambda)>\omega$ were obtained in \cite{MR0376359}, shortly before the discovery of pcf theory by Shelah. However, the methods of \cite{MR0376359} require uncountable cofinality.} namely $|pcf(\mathfrak{a})|\leq 2^{|\mathfrak{a}|}$.
But this bound is not absolute, as $2^{|\mathfrak{a}|}$ can be easily manipulated by forcing.
Thus the real interesting bound is the $\aleph$-scale bound, which says that $|pcf(\mathfrak{a})|<|\mathfrak{a}|^{+4}$.
We indicate, however, that sometimes (e.g., if $2^{|\mathfrak{a}|}=|\mathfrak{a}|^+$) the $\beth$-scale bound is better than the $\aleph$-scale bound.

The proof of the $\aleph$-scale bound is a combination of basic properties of $pcf(\mathfrak{a})$ and a prediction principle called \emph{club guessing}.
The latter will be discussed extensively in the next section, so here we conclude with one important feature of $pcf(\mathfrak{a})$, called \emph{localization}.
Needless to say that this property of $pcf(\mathfrak{a})$ was proved by Shelah.

\begin{theorem}
  \label{thmlocalization} Let $\mathfrak{a}$ be progressive and assume that $\mathfrak{b}\subseteq pcf(\mathfrak{a})$ is also progressive.
  Let $\lambda\in pcf(\mathfrak{b})$.
  Then one can find $\mathfrak{b}_0\subseteq\mathfrak{b}$ such that $|\mathfrak{b}_0|\leq|\mathfrak{a}|$ and $\lambda\in pcf(\mathfrak{b}_0)$.
\end{theorem}

Intuitively, this property puts a restriction on the size of $pcf(\mathfrak{a})$.
It is closely related to the Achilles and the Tortoise property, which says that within $pcf(\mathfrak{a})$ one cannot construct long subsets of the from $\mathfrak{b}$ for which $\max pcf(\mathfrak{b}\cap\lambda)>\lambda$.
The bound of $|pcf(\mathfrak{a})|<|\mathfrak{a}|^{+4}$ can be derived from the localization property in the context of an appropriate instance of club guessing, as will be shown later.
For more background in pcf theory we suggest \cite{MR2768693} and \cite{MR1086455}, as well as Shelah's monograph \cite{MR1318912}.

\newpage

\section{Club guessing and Galvin's property}

Let $\kappa,\lambda$ be regular cardinals\footnote{In the general setting, $\lambda$ can be singular provided that $\cf(\lambda)\geq\kappa^{++}$. But for the applications to cardinal arithmetic we may restrict our attention to $\lambda=\cf(\lambda)$.} where $\kappa<\lambda$.
The set $S^\lambda_\kappa$ is a stationary subset of $\lambda$.
Fix $S\subseteq S^\lambda_\kappa$ that is stationary in $\lambda$, and let $\mathcal{C}=(C_\delta\mid\delta\in{S})$ be a sequence of sets such that $C_\delta\subseteq\delta$ is a club of $\delta$ of order-type $\kappa$ for every $\delta\in{S}$.
One says that $\mathcal{C}$ is a club guessing sequence for $\lambda$ iff for every club $E$ of $\lambda$ one can find $\delta\in{S}$ for which $C_\delta\subseteq{E}$.
The following theorem appears in \cite{MR1318912}, and it plays a key role in proofs which bound the cardinality of $pcf(\mathfrak{a})$.

\begin{theorem}
  \label{thmcgsaharon} Suppose that $\lambda\geq\cf(\lambda)\geq\kappa^{++}$, where $\kappa$ is a regular cardinal.
  Let $S\subseteq S^\lambda_\kappa$ be a stationary subset of $\lambda$.
  Then there exists a club guessing sequence $(C_\alpha\mid\alpha\in{S})$.
\end{theorem}

It is important to notice that there is a gap here between the small parameter $\kappa$ (that is, the size of each $C_\alpha$) and the cofinality of $\lambda$ (that is, $\kappa^{++}$).
In the typical case of $\lambda=\kappa^{++}$ this gap reduces to two cardinalities, and in general this optimal gap cannot be avoided.
Indeed, one can force the failure of club guessing at $S^{\kappa^+}_\kappa$.
As we shall see later, this gap is instrumental when trying to compute the size of $pcf(\mathfrak{a})$.
Our goal in this section is to obtain club guessing at $S^{\kappa^+}_\kappa$ from an instance of Galvin's property.

Galvin showed\footnote{The proof was published in \cite{MR0369081}.} that if $\kappa=\kappa^{<\kappa}>\aleph_0$ and $\mathscr{F}$ is a normal filter over $\kappa$ then every family $\mathcal{C}=\{C_\gamma\mid\gamma\in\kappa^+\}\subseteq\mathscr{F}$ admits a subfamily $\{C_{\gamma_i}\mid i\in\kappa\}$ such that $\bigcap\{C_{\gamma_i}\mid i\in\kappa\}\in\mathscr{F}$.
A central example is the club filter over $\kappa$, denoted by $\mathscr{D}_\kappa$.
We use the notation ${\rm Gal}(\mathscr{F},\kappa,\kappa^+)$ to denote the above statement.
The assumption $\kappa=\kappa^{<\kappa}$ boils down to a local instance of \textsf{GCH} if one wishes to apply Galvin's theorem at successor cardinals.
Thus if $2^\lambda=\lambda^+$ then ${\rm Gal}(\mathscr{D}_{\lambda^+},\lambda^+,\lambda^{++})$ holds true.

A natural question is whether Galvin's assumption $2^\lambda=\lambda^+$ is droppable.
The answer turns out to be interesting.
Abraham and Shelah proved in \cite{MR830084} that Galvin's property consistently fails.
More specifically, they proved that if $\kappa$ is regular and $\lambda\geq\cf(\lambda)>\kappa^+$ then one can force $2^{\kappa^+}=\lambda$ with $\neg{\rm Gal}(\mathscr{D}_{\kappa^+},\kappa^+,\lambda)$.
This statement is called \emph{the ultimate failure} of Galvin's property at $\mathscr{D}_{\kappa^+}$, since the size of the family witnessing the failure of Galvin's property is the largest possible (namely, $2^{\kappa^+}$).
We conclude, therefore, that Galvin's theorem is not a \textsf{ZFC} statement.

However, one can prove instances of Galvin's property under weaker assumptions than the one used by Galvin.
It was shown in \cite{MR3604115} that if $2^\kappa<2^{\kappa^+}=\lambda$ then ${\rm Gal}(\mathscr{D}_{\kappa^+},\kappa^+,\lambda)$ holds.\footnote{Recall that the statement $2^\kappa<2^{\kappa^+}$ is equivalent to the prediction principle $\Phi_{\kappa^+}$.}
Thus $2^\kappa<2^{\kappa^+}$ yields instances of Galvin's property even if $2^\kappa>\kappa^+$, and the original assumption of Galvin can be relaxed.
The crucial point here becomes meaningful if one adopts a global point of view.
Galvin's assumption $2^\kappa=\kappa^+$ may fail everywhere, while $\Phi_\kappa$ must hold at many places, as shown in the previous section.
This point will be relevant in the sequel, due to the main result of this section which derives club guessing (with a gap of one cardinal) from an instance of Galvin's property.

\begin{theorem}
  \label{thmcggalvin} Let $\kappa$ be a regular and uncountable cardinal, and let $S\subseteq S^{\kappa^+}_\kappa$ be stationary.
  If ${\rm Gal}(\mathscr{D}_{\kappa^+},\kappa^+,2^{\kappa^+})$ holds then there exists a club guessing sequence $(C_\alpha\mid\alpha\in{S})$.
\end{theorem}

\par\noindent\emph{Proof}. \newline
Let $\mathcal{C}=(C_\delta\mid\delta\in{S})$ be a sequence of sets, where $C_\delta\subseteq\delta$ is a club of $\delta$ for each $\delta\in{S}$.
Let $E\subseteq\kappa^+$ be a club of $\kappa^+$.
Denote by $\mathcal{C}\upharpoonright{E}$ the sequence $(C_\delta\cap{E}\mid\delta\in acc(E)\cap{S})$.
We move from $E$ to $acc(E)$ in order to make sure that $C_\delta\cap{E}$ is a club of $\delta$.

We claim that for some $E\subseteq\kappa^+$, the sequence $\mathcal{C}\upharpoonright{E}$ is a club guessing sequence.\footnote{The index set of $\mathcal{C}\upharpoonright{E}$ is $acc(E)\cap{S}\subseteq{S}$, but trivially if there is a club guessing sequence for a stationary subset of $S$ then there is such a sequence for $S$.}
Assume towards contradiction that this claim fails.
Therefore, for every club $E$ of $\kappa^+$ there exists a club $D_E\subseteq\kappa^+$ such that if $\delta\in acc(E)\cap{S}$ then $C_\delta\cap{E}\nsubseteq D_E$.
Of course, if $D\subseteq{D_E}$ then $C_\delta\cap{E}\nsubseteq{D}$, so we may shrink each $D_E$ by letting $A_E=E\cap{D_E}$, and now for every $\delta\in acc(E)\cap{S}$ one has $C_\delta\cap{E}\nsubseteq A_E$.

Let $\mathcal{A}=\{A_E\mid E\ \text{is a club of}\ \kappa^+\}$.
We may assume, without loss of generality,\footnote{Our goal is to apply ${\rm Gal}(\mathscr{D}_{\kappa^+},\kappa^+,2^{\kappa^+})$ to $\mathcal{A}$. Now if $|\mathcal{A}|<2^{\kappa^+}$ then there are $2^{\kappa^+}$ many $E$s with the same $A_E$, and they satisfy (trivially) the conclusion of the Galvin property, so we may assume that $|\mathcal{A}|=2^{\kappa^+}$.} that $|\mathcal{A}|=2^{\kappa^+}$.
Notice that each element of $\mathcal{A}$ is in $\mathscr{D}_{\kappa^+}$, hence ${\rm Gal}(\mathscr{D}_{\kappa^+},\kappa^+,2^{\kappa^+})$ applies.
Fix a family $\mathcal{A}'=\{A_{E_i}\mid i\in\kappa^+\}\subseteq\mathcal{A}$ and a club $E_*\subseteq\kappa^+$ such that $E_*\subseteq A_{E_i}$ for every $i\in\kappa^+$.
If $\delta\in acc(E_*)\cap{S}$ then $C_\delta\cap E_*$ is a club of $\delta$.
Since $C_\delta\cap{E_*}\subseteq C_\delta\cap{A_{E_i}}$ we conclude that $C_\delta\cap{A_{E_i}}$ is a club of $\delta$ as well, for every $i\in\kappa^+$.
In particular, $|C_\delta\cap{A_{E_i}}|=\kappa$.
By induction on $\alpha\in\kappa^+$ we define $E^\alpha\subseteq\kappa^+$ as follows.
We let $E^0=A_{E_0}$ and $E^\alpha=\bigcap_{\beta\in\alpha}E^\beta$ whenever $\alpha$ is a limit ordinal of $\kappa^+$.
Finally, $E^{\alpha+1}=E^\alpha\cap A_{E_\alpha}$ for every $\alpha\in\kappa^+$.
Notice that each $E^\alpha$ is a club of $\kappa^+$.

Fix $\delta\in acc(E_*)\cap{S}$.
Since $E_*\subseteq A_{E_\alpha}$ for every $\alpha\in\kappa^+$ one concludes that $C_\delta\cap{E_*}\subseteq C_\delta\cap{A_{E_\alpha}}$ for every $\alpha\in\kappa^+$.
The sequence $(C_\delta\cap{E^\alpha}\mid\alpha\in\kappa^+)$ is $\subseteq$-decreasing, and since $|C_\delta\cap{E^\alpha}|=\kappa$ for every $\alpha\in\kappa^+$ we see that there is some $\alpha\in\kappa^+$ for which $C_\delta\cap{E^\alpha}=C_\delta\cap{E_*}$.
It follows that $C_\delta\cap{E^\alpha}=C_\delta\cap{E^{\alpha+1}}$ since $E^\alpha\supseteq E^{\alpha+1}\supseteq E_*$, so  $C_\delta\cap{E^\alpha}=C_\delta\cap{E^{\alpha+1}}\subseteq E^{\alpha+1}$.
On the other hand, $E^{\alpha+1}\subseteq A_{E_\alpha}$ and $\delta\in acc(E^\alpha)$ so $C_\delta\cap{E^\alpha}\nsubseteq E^{\alpha+1}$, a contradiction.

\hfill \qedref{thmcggalvin}

We conclude this section with a comparison between the classical club guessing of Shelah (i.e., Theorem \ref{thmcgsaharon}) and the current version which comes from the Galvin property.
In both cases we have two parameters.
The smaller parameter $\kappa$ is regular and if forms the size of the guessing club $C_\delta$.
The bigger parameter $\lambda$ is the domain of the clubs to be guessed.
In both proofs one begins with an arbitrary sequence $\mathcal{C}=(C_\delta\mid\delta\in{S})$, where $S\subseteq S^\lambda_\kappa$ is stationary and each $C_\delta$ is a club of $\delta$.
Then one argues that $\mathcal{C}\upharpoonright{E}$ is a club guessing sequence for some $E\subseteq\lambda$.

In order to define $E$ one creates a decreasing sequence $(E^\alpha\mid\alpha\in\kappa^{+}+1)$ of clubs of $\lambda$, where the last element is $E=E^{\kappa^+}$.
The sequence is continuous and hence $E$ is the intersection of $\kappa^+$-many clubs of $\lambda$.
Apart from the last step of $E$, the length of the sequence is $\kappa^+$, and this is necessary in order to stabilize the decreasing derived sequences $(C_\delta\cap{E^\alpha}\mid\alpha\in\kappa^+)$, bearing in mind that $otp(C_\delta)=\kappa$.
The last step in which $E=E^{\kappa^+}$ is created, compels $\lambda$ to be $\kappa^{++}$ (or, more generally, $\cf(\lambda)\geq\kappa^{++}$).
This is the reason for the gap between $\kappa$ and $\lambda$ in the classical club guessing theorem.

But in the presence of Galvin's property one obtains a decreasing sequence of clubs of $\lambda$ so that every element in the sequence contains a fixed club $E_*$, which serves (at the end) as $E^{\kappa^+}$.
This can be done, under a mild assumption, even for $\lambda=\kappa^+$.
Put another way, Galvin's property helps to reduce the gap between $\kappa$ and $\lambda$ to one cardinality.
Thus one obtains club guessing at $S^{\kappa^+}_\kappa$, and this is crucial for computing the size of $pcf(\mathfrak{a})$ as we shall see in the next section.

\newpage

\section{A legend of three and four}

There are three things which are stately in their march, four which are stately in their going.\footnote{See \cite[Chapter 30:29]{proverbs}.}
Our goal in this section is to supply a mathematical interpretation to the above quotation.
Familiarity with the proof of Shelah's bound $|pcf(\mathfrak{a})|<|\mathfrak{a}|^{+4}$ where $\mathfrak{a}$ is a progressive interval of regular cardinals leads to the conclusion that an improvement in club guessing would give a better bound on the size of $pcf(\mathfrak{a})$.
In particular, if $|\mathfrak{a}|=\eta$ then club guessing at $S^{\eta^{++}}_{\eta^+}$ yields the corresponding bound of $|pcf(\mathfrak{a})|<|\mathfrak{a}|^{+3}$.

In this section we spell-out the proof of this statement.
The only deviation from the classical proof of Shelah is when we replace the \textsf{ZFC} club guessing by the stronger version based on the Galvin property.
Hence we will be able to improve the bound on the size of $pcf(\mathfrak{a})$ once we show that an appropriate instance of Galvin's property holds at relevant places, and this will be done later.
The proof of the following theorem is based on the presentation in \cite{MR1086455}.
We indicate that one can use the ideas of \cite{MR2768693} as well.

\begin{theorem}
  \label{thmpcfbound3} Let $\mathfrak{a}$ be a progressive interval of regular cardinals, and let $\eta=|\mathfrak{a}|$.
  Suppose that there is a club guessing sequence at $S^{\eta^{++}}_{\eta^+}$.
  Then $|pcf(\mathfrak{a})|<|\mathfrak{a}|^{+3}$.
\end{theorem}

\par\noindent\emph{Proof}. \newline
By omitting one element from $\mathfrak{a}$ (if needed) we may assume that $\min(\mathfrak{a})$ is a successor cardinal, say $\aleph_{\delta+1}$.
It is easy to see that (under this slight modification of $\mathfrak{a}$) all the elements of $pcf(\mathfrak{a})$ are successor cardinals.
Denote $\max pcf(\mathfrak{a})$ by $\aleph_{\delta+\rho+1}$.

We shall define a topological structure whose underlying set is $\rho+1$, and then we will show that the properties of this topological space imply that $|\rho|\leq\eta^{++}$.
In order to generate our topology we define the following closure operation.
Given $x\subseteq\rho+1$ we let:
\[c\ell(x)=\{\gamma\leq\rho\mid\aleph_{\delta+\gamma+1}\in pcf(\{\aleph_{\delta+\beta+1}\mid\beta\in{x}\})\}.\]
One can verify that $c\ell(\varnothing)=\varnothing$, that $x\subseteq c\ell(x)$ for every $x$, that $x\subseteq{y}$ implies $c\ell(x)\subseteq c\ell(y)$, that $c\ell(x\cup{y})=c\ell(x)\cup c\ell(y)$ and that $c\ell(c\ell(x))=c\ell(x)$.
Less routine properties come from the attributes of pcf, in particular:
\begin{enumerate}
  \item [$(a)$] If $x\subseteq\rho+1$ and $\gamma\in c\ell(x)$ then there exists $y\subseteq{x}$ such that $|y|\leq\eta$ and $\gamma\in c\ell(y)$.
  \item [$(b)$] For every $x\subseteq\rho+1$ there is a last element in $c\ell(x)$.
  \item [$(c)$] If $\omega<\cf(\gamma)\leq\gamma\leq\rho$ then there is a club $c\subseteq\gamma$ such that $c\ell(c)\subseteq\gamma+1$.
\end{enumerate}

Observe that $(a)$ is simply the localization property, and $(b)$ is the fact that $pcf(\mathfrak{b})$ has a last element.

Assume towards contradiction that $|\rho|\geq\eta^{+3}$, so without loss of generality $|\rho|=\eta^{+3}$.
Let $S=S^{\eta^{++}}_{\eta^+}$ and let $(C_\alpha\mid\alpha\in{S})$ be a club guessing sequence.
Fix a sufficiently large regular cardinal $\chi$.
Let $(M_\beta\mid\beta\leq\eta^{++})$ be an increasing continuous sequence of elementary submodels of $\mathcal{H}(\chi)$ for which the following requirements are met:
\begin{enumerate}
  \item [$(\aleph)$] $\eta^{++}\subseteq M_\beta$ and $|M_\beta|=\eta^{++}$ for each $\beta\in\eta^{++}$.
  \item [$(\beth)$] $(M_\gamma\mid\gamma\leq\beta)\in M_{\beta+1}$ for every $\beta\in\eta^{++}$.
  \item [$(\gimel)$] $(C_\alpha\mid\alpha\in{S})\in M_\beta$ for every $\beta\in\eta^{++}$.
  \item [$(\daleth)$] $\{\langle x,c\ell(x)\rangle\mid x\subseteq\eta^{+3}+1\}\in M_\beta$ for every $\beta\in\eta^{++}$.
\end{enumerate}
For each $\beta\in\eta^{++}+1$ let $\gamma_\beta=M_\beta\cap\eta^{+3}$, so $\gamma_\beta\in\eta^{+3}$.
It follows that $D=\{\gamma_\beta\mid\beta\leq\eta^{++}\}$ is a closed bounded subset of $\eta^{+3}$.
Notice that $(\gamma_\delta\mid\delta\in\beta)\in M_{\beta+1}$ for every $\beta\in\eta^{++}$ since each $\gamma_\delta$ belongs to $M_{\beta+1}$ (being definable in $M_{\beta+1}$) and then the whole sequence $(\gamma_\delta\mid\delta\in\beta)$ is in $M_{\beta+1}$ by virtue of $(\beth)$.

For every $\alpha,\beta\in\eta^{++}$ let $E^\beta_\alpha=\{\gamma_\delta\mid\delta\in C_\alpha\cap\beta\}$.
Observe that $E^\beta_\alpha$ is definable in $M_{\beta+1}$, and hence $E^\beta_\alpha\in M_{\beta+1}$.
Therefore, if $c\ell(E^\beta_\alpha)$ is bounded in $\eta^{+3}$ then this bound is computable in $M_{\beta+1}$ and hence belongs to $M_{\beta+1}$.
It follows that if $c\ell(E^\beta_\alpha)$ is bounded in $\eta^{+3}$ then $c\ell(E^\beta_\alpha)\subseteq\gamma_{\beta+1}$.

From property $(c)$ one infers that there exists a club $E$ of $\gamma_{\eta^{++}}$ such that $c\ell(E)\subseteq\gamma_{\eta^{++}}+1$.
For every $\alpha\in{S}$ let $T_\alpha=\{\gamma_\beta\mid\beta\in{C_\alpha}\}$, so $T_\alpha$ is the $D$-copy of the element $C_\alpha$ in the club guessing sequence.
Since $C_\alpha$ is a club of $\alpha$, $T_\alpha$ is a club of $\gamma_\alpha$.
Fix an ordinal $\alpha$ so that $T_\alpha\subseteq{E}$.
Let $\zeta$ be the last element of $c\ell(T_\alpha)$, it exists by $(b)$.
Observe that $\zeta\geq\gamma_\beta$ for every $\gamma_\beta\in T_\alpha$ and hence $\zeta\geq\gamma_\alpha$.

On the other hand, there must be some $\beta\in\alpha$ for which $\zeta\in c\ell(T_\alpha\cap\gamma_\beta)$.
Indeed, $|T_\alpha|=\eta^+$.
Hence, if $y\subseteq T_\alpha, |y|=\eta$ and $\zeta\in c\ell(y)$ then $y\cap T_\alpha$ is bounded in $\gamma_\alpha$ (recall that $\cf(\gamma_\alpha)=\eta^+$), thus $y=y\cap T_\alpha$ is a subset of $T_\alpha\cap\gamma_\beta$ for some $\beta\in\alpha$.
Therefore, one concludes that $\zeta\in c\ell(T_\alpha\cap\gamma_\beta)$ for some $\beta\in\alpha$, as wanted.
Using the above notation, $T_\alpha\cap\gamma_\beta=E^\beta_\alpha$ and hence $c\ell(E^\beta_\alpha)\subseteq c\ell(T_\alpha)\subseteq c\ell(E)\subseteq\gamma_{\eta^{++}}+1$.
In particular, $c\ell(E^\beta_\alpha)$ is bounded in $\eta^{+3}$.
As indicated above, the bound belongs to $M_{\beta+1}$.
Therefore, $\zeta\in\gamma_{\beta+1}<\gamma_\alpha$, a contradiction.

\hfill \qedref{thmpcfbound3}

Here is an easy conclusion which gives an improved pcf bound:

\begin{corollary}
  \label{corweakdiamond} Assume that:
  \begin{enumerate}
    \item [$(a)$] $\lambda=\aleph_\delta$ is a strong limit singular cardinal.
    \item [$(b)$] $\lambda$ is not a fixed point of the $\aleph$-function.
    \item [$(c)$] $\mathfrak{a}$ is an end-segment of ${\rm Reg}\cap\lambda, |\mathfrak{a}|=\eta$.
    \item [$(d)$] $2^{\eta^+}<2^{\eta^{++}}$.
  \end{enumerate}
  Then $|pcf(\mathfrak{a})|<|\mathfrak{a}|^{+3}$, and hence $2^\lambda<\aleph_{|\delta|^{+3}}$.
\end{corollary}

\par\noindent\emph{Proof}. \newline
By \cite{MR3604115}, assumption $(d)$ implies ${\rm Gal}(\mathscr{D}_{\eta^{++}},\eta^{++},2^{\eta^{++}})$.
Since $\eta^+$ is regular we conclude from Theorem \ref{thmcggalvin} that there exists a club guessing sequence at $S^{\eta^{++}}_{\eta^+}$.
Applying Theorem \ref{thmpcfbound3} we see that $|pcf(\mathfrak{a})|<|\mathfrak{a}|^{+3}$.
Since $\lambda$ is a strong limit cardinal, $2^\lambda=\max pcf(\mathfrak{a})$.
Therefore, $2^\lambda<\aleph_{|\delta|^{+3}}$ as desired.

\hfill \qedref{corweakdiamond}

The corollary shows that an instance of weak diamond yields a locally interesting pcf bound.
But the real import of the weak diamond hinges upon the fact that it holds at unboundedly many points in \textsf{ZFC}, as shown before.
Of course, we need instances of weak diamond at double successors, and the existence of these instances is not a \textsf{ZFC} statement.
However, mild assumptions produce the desired setting.

\begin{theorem}
  \label{thmglobalpcfbound3} Assume that every limit cardinal is a strong limit cardinal.
  Then there is a class $\mathfrak{C}$ of singular cardinals such that:
  \begin{enumerate}
    \item [$(\aleph)$] If $\aleph_\delta\in\mathfrak{C}$ and $\mathfrak{a}\subseteq{\rm Reg}\cap\aleph_\delta$ is progressive then $|pcf(\mathfrak{a})|<|\mathfrak{a}|^{+3}$.
    \item [$(\beth)$] If $\aleph_\delta\in\mathfrak{C}$ then $2^{\aleph_\delta}<\aleph_{(|\delta|^{+3})}$.
  \end{enumerate}
\end{theorem}

\par\noindent\emph{Proof}. \newline
Fix an infinite cardinal $\kappa_0$ and let $\lambda_0$ be the first singular cardinal greater than $\kappa_0$.
Let $\kappa=\kappa_0^+$.
Since $\lambda_0$ is strong limit, $\mu=2^\kappa<\lambda_0$.
By the arguments of Theorem \ref{thmclassofwd} there must be some $\chi\geq\kappa^+$ such that $\chi\leq\mu$ and $\Phi_\chi$ holds.
Notice that $\chi$ is necessarily a double successor cardinal, so from Corollary \ref{corweakdiamond} one can find a singular strong limit cardinal $\aleph_\delta$ such that $\delta<\aleph_\delta$ and a progressive interval $\mathfrak{a}\subseteq{\rm Reg}\cap\aleph_\delta$ so that $|\mathfrak{a}|=\chi$, and infer that $|pcf(\mathfrak{a})|<|\mathfrak{a}|^{+3}$.
Since $\aleph_\delta$ is a strong limit cardinal, $2^{\aleph_\delta}=\max pcf(\mathfrak{a})$ and hence $2^{\aleph_\delta}<\aleph_{(|\delta|^{+3})}$.
This reasoning holds with respect to every starting point $\kappa_0$, so we are done.

\hfill \qedref{thmglobalpcfbound3}

The assumption that every limit cardinal is a strong limit cardinal is much stronger than the assumption needed for getting the conclusion of the above theorem.
Basically, in order to avoid instances of $2^{\eta^+}<2^{\eta^{++}}$ one has to accept a restricted constellation of cardinal arithmetic.
Let us describe a typical cardinal arithmetic setting in which Corollary \ref{corweakdiamond} does not apply.
In this setting, all the relevant instances of weak diamond concentrate on weakly but not strongly inaccessible cardinals.\footnote{Of course, there might be weak diamonds on large cardinals, but these are irrelevant to our arguments, since by starting from any $\kappa_0$ we consider possible instances of weak diamond between $\kappa_0$ and $2^{\kappa_0}$.}
In particular, there is a class of weakly but not strongly inaccessible cardinals and $2^{\eta^+}$ is weakly inaccessible for every successor cardinal of the form $\eta^+$.
Moreover, $2^{\eta^+}=2^{\eta^{++}}$ everywhere, so cardinal arithmetic is very peculiar: there are long intervals of regular cardinals with a constant value of the power set, and all the regular cardinals are arranged in this way.

This is not the only setting in which Corollary \ref{corweakdiamond} may fail, but the other possibilities are similar.
That is, long intervals with the same power set, with values either at weakly but not strongly inaccessible cardinals, or at their successors, or at successors of singular cardinals.
The main thing is that the assumption that $2^{\eta^+}=2^{\eta^{++}}$ holds everywhere is quite restricting.
We state, therefore, the following:

\begin{conjecture}
  \label{conjbound3} In any model of \textsf{ZFC} there is a class of strong limit singular cardinals of the form $\aleph_\delta$ for which $2^{\aleph_\delta}<\aleph_{(|\delta|^{+3})}$.
\end{conjecture}

\newpage

\bibliographystyle{alpha}
\bibliography{arlist}

\end{document}